\documentclass[a4paper,11pt]{article}

\usepackage{graphicx}
\usepackage{amssymb,amsmath,amsthm,times,mathrsfs,fullpage}
\usepackage{tikz,wrapfig,multirow}
\usepackage[pdftex]{hyperref}
\hypersetup{plainpages=True, pdfstartview=FitH, bookmarksopen=true, colorlinks=true,linkcolor=blue,citecolor=blue}

\def\pf{\noindent \emph{Proof.}\ }
\def\qed{{\quad\rule{1mm}{3mm}\,}}

\begin{document}

\pgfdeclarelayer{background}
\pgfdeclarelayer{foreground}
\pgfsetlayers{background,main,foreground}

\newtheorem{thm}{Theorem}
\newtheorem{cor}{Corollary}
\newtheorem{lmm}{Lemma}
\newtheorem{conj}{Conjecture}
\newtheorem{pro}{Proposition}
\newtheorem{Def}{Definition}
\theoremstyle{remark}\newtheorem{Rem}{Remark}

\title{On the Asymptotic Growth of the Number of Tree-Child Networks}
\author{Michael Fuchs\thanks{MF was partially supported by grant MOST-107-2115-M-009-
010-MY2.}\\
    Department of Mathematical Sciences\\
    National Chengchi University\\
    Taipei, 116\\
    Taiwan \and
    Guan-Ru Yu\thanks{GRY was partially supported by FWF SFB F50-10.}\\
    Faculty of Mathematics\\
    University of Vienna\\
    Vienna, 1090\\
    Austria \and
    Louxin Zhang\thanks{LZ was financially supported by Singapore MOE Tier-1 Research Fund.}\\
    Department of Mathematics\\
    National University of Singapore\\
    Singapore, 119076\\
    Singapore}

\maketitle

\begin{abstract}
In a recent paper, McDiarmid, Semple, and Welsh (2015) showed that the number of tree-child networks with $n$ leaves has the factor $n^{2n}$ in its main asymptotic growth term. In this paper, we improve this by completely identifying the main asymptotic growth term up to a constant. More precisely, we show that the number of tree-child networks with $n$ leaves grows like
\[
\Theta\left(n^{-2/3}e^{a_1(3n)^{1/3}}\left(\frac{12}{e^2}\right)^{n}n^{2n}\right),
\]
where $a_1=-2.338107410\cdots$ is the largest root of the Airy function of the first kind. For the proof, we bijectively map the underlying graph-theoretical problem onto a problem on words. For the latter, we can find a recurrence to which a recent powerful asymptotic method of Elvey Price, Fang, and Wallner (2019) can be applied.
\end{abstract}

\section{Introduction}

Over the last two decades, {\it phylogenetic networks} have become increasingly popular and have been used more and more frequently in modeling horizontal genetic transfer events in evolutionary genomics. Because of their now widespread usage, studying basic combinatorial properties such as counting them has attracted some recent efforts; see, e.g., Bouvel et al. \cite{BoGaMa}, Cardona and Zhang \cite{Zh2}, Fuchs et al. \cite{FuGiMa}, Gunawan et al. \cite{Gunawan_20},  McDiarmid at al. \cite{DiSeWe}, and Zhang \cite{Zh1}. While the counting of phylogenetic trees goes back at least to Schr\"{o}der \cite{Sch} and we have a complete combinatorial understanding of it, still very little is known about combinatorial counting questions for phylogenetic networks and their subclasses.

A (rooted, binary, leaf-labeled) phylogenetic network with $n$ leaves is defined as a rooted directed acyclic graph (or DAG for short) which is connected and has no parallel edges and whose vertices fall into one of the following four categories:
\begin{itemize}
\item[(i)] a root $\rho$ of indegree $0$ and outdegree $1$;
\item[(ii)] $n$ leaves of indegree $1$ and outdegree $0$ which are bijectively labeled by elements from the set $\{1,\ldots, n\}$;
\item[(iii)] nodes of indegree $1$ and outdegree $2$ which are called {\it tree nodes};
\item[(iv)] nodes of indegree $2$ and outdegree $1$ which are called {\it reticulation nodes}.
\end{itemize}

Many subclasses of phylogenetic networks have been considered. An important one arising from phylogenetic applications is the class containing {\it tree-child networks} which we will define next.
\begin{Def}\label{def-tc}
A phylogenetic network is called a tree-child network if every node which is not a leaf has at least one child which is not a reticulation node.
\end{Def}

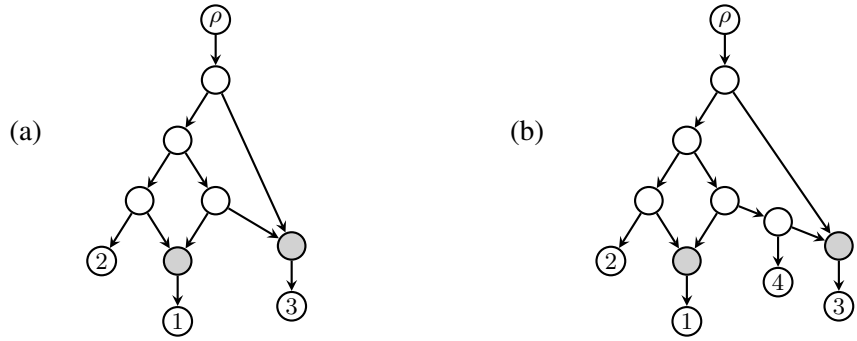
\begin{figure}[h]
\begin{center}
\begin{tikzpicture}
\draw (-2.5cm,-1.5cm) node (0) {(a)};
\draw (0cm,0cm) node[inner sep=1.2pt,line width=0.8pt,draw,circle] (1) {{\footnotesize $\rho$}};
\draw (0cm,-0.8cm) node[minimum size=7.8pt,line width=0.8pt,draw,circle] (2) {};
\draw (-0.5cm,-1.6cm) node[minimum size=7.8pt,line width=0.8pt,draw,circle] (3) {};
\draw (1cm,-3cm) node[minimum size=7.8pt,line width=0.8pt,draw,circle,fill=gray!35] (6) {};
\draw (-1cm,-2.4cm) node[minimum size=7.8pt,line width=0.8pt,draw,circle] (4) {};
\draw (0cm,-2.4cm) node[minimum size=7.8pt,line width=0.8pt,draw,circle] (5) {};
\draw (-1.5cm,-3.2cm) node[inner sep=1.2pt,line width=0.8pt,draw,circle] (7) {{\footnotesize $2$}};
\draw (-0.5cm,-3.2cm) node[minimum size=7.8pt,line width=0.8pt,draw,circle,fill=gray!35] (8) {};
\draw (-0.5cm,-4cm) node[inner sep=1.2pt,line width=0.8pt,draw,circle] (9) {{\footnotesize $1$}};
\draw (1cm,-3.8cm) node[inner sep=1.2pt,line width=0.8pt,draw,circle] (10) {{\footnotesize $3$}};

\draw[-stealth,line width=0.8pt] (1) -- (2);
\draw[-stealth,line width=0.8pt] (2) -- (3);
\draw[-stealth,line width=0.8pt] (2) -- (6);
\draw[-stealth,line width=0.8pt] (3) -- (4);
\draw[-stealth,line width=0.8pt] (3) -- (5);
\draw[-stealth,line width=0.8pt] (4) -- (8);
\draw[-stealth,line width=0.8pt] (5) -- (8);
\draw[-stealth,line width=0.8pt] (5) -- (6);
\draw[-stealth,line width=0.8pt] (4) -- (7);
\draw[-stealth,line width=0.8pt] (8) -- (9);
\draw[-stealth,line width=0.8pt] (6) -- (10);

\draw (4.1cm,-1.5cm) node (0) {(b)};
\draw (6.7cm,0cm) node[inner sep=1.2pt,line width=0.8pt,draw,circle] (1) {{\footnotesize $\rho$}};
\draw (6.7cm,-0.8cm) node[minimum size=7.8pt,line width=0.8pt,draw,circle] (2) {};
\draw (6.2cm,-1.6cm) node[minimum size=7.8pt,line width=0.8pt,draw,circle] (3) {};
\draw (8.2cm,-3cm) node[minimum size=7.8pt,line width=0.8pt,draw,circle,fill=gray!35] (6) {};
\draw (7.4cm,-2.68cm) node[minimum size=7.8pt,line width=0.8pt,draw,circle] (11) {};
\draw (7.4cm,-3.48cm) node[inner sep=1.2pt,line width=0.8pt,draw,circle] (12) {{\footnotesize $4$}};
\draw (5.7cm,-2.4cm) node[minimum size=7.8pt,line width=0.8pt,draw,circle] (4) {};
\draw (6.7cm,-2.4cm) node[minimum size=7.8pt,line width=0.8pt,draw,circle] (5) {};
\draw (5.2cm,-3.2cm) node[inner sep=1.2pt,line width=0.8pt,draw,circle] (7) {{\footnotesize $2$}};
\draw (6.2cm,-3.2cm) node[minimum size=7.8pt,line width=0.8pt,draw,circle,fill=gray!35] (8) {};
\draw (6.2cm,-4cm) node[inner sep=1.2pt,line width=0.8pt,draw,circle] (9) {{\footnotesize $1$}};
\draw (8.2cm,-3.8cm) node[inner sep=1.2pt,line width=0.8pt,draw,circle] (10) {{\footnotesize $3$}};

\draw[-stealth,line width=0.8pt] (1) -- (2);
\draw[-stealth,line width=0.8pt] (2) -- (3);
\draw[-stealth,line width=0.8pt] (2) -- (6);
\draw[-stealth,line width=0.8pt] (3) -- (4);
\draw[-stealth,line width=0.8pt] (3) -- (5);
\draw[-stealth,line width=0.8pt] (4) -- (8);
\draw[-stealth,line width=0.8pt] (5) -- (8);
\draw[-stealth,line width=0.8pt] (5) -- (11);
\draw[-stealth,line width=0.8pt] (11) -- (6);
\draw[-stealth,line width=0.8pt] (11) -- (12);
\draw[-stealth,line width=0.8pt] (4) -- (7);
\draw[-stealth,line width=0.8pt] (8) -- (9);
\draw[-stealth,line width=0.8pt] (6) -- (10);
\end{tikzpicture}
\end{center}
\caption{(a) A phylogenetic network which is not a tree-child network because its two reticulations nodes (in gray) are the children of the same tree node; (b) A tree-child network.}\label{pt-and-tc}
\end{figure}

See Figure~\ref{pt-and-tc} for an example of a phylogenetic network which is not a tree-child network (a) and a phylogenetic network which is a tree-child network (b). 

We denote the set of tree-child networks with $n$ leaves by $\mathcal{TC}_n$ and its cardinality by ${\rm TC}_n$ throughout this work. The first few terms of the latter sequence were computed, e.g., in \cite{Zh2}:
\[
\{{\rm TC}_n\}_{n\geq 2}=\{3, 66, 4059, 496710, 101833875, 31538916360, \ldots\}.
\]

\vspace*{0.2cm}

In order to understand the growth of this sequence, asymptotic counting results have been proved. The first and in fact still best such result was obtained in \cite{DiSeWe}, where the authors showed that there exist constants $0<c_1<c_2$ such that for all $n$:
\[
(c_1n)^{2n}\leq{\rm TC}_{n}\leq (c_2n)^{2n}.
\]
This results identifies $n^{2n}$ as the main term in the asymptotic growth of ${\rm TC}_n$. However, it does not yield anything even for the exponential growth rate since no explicit values for $c_1,c_2$ were given in \cite{DiSeWe}.

It is the purpose of this paper to improve upon this result. More precisely, we will find all terms in the main asymptotic growth term of the above counting sequence.

\begin{thm}\label{main-result}
The number of tree-child networks with $n$ leaves satisfies
\[
{\rm TC}_{n}=\Theta\left(n^{-2/3}e^{a_1(3n)^{1/3}}\left(\frac{12}{e^2}\right)^{n}n^{2n}\right),
\]
where $a_1$ is the largest root of the Airy function ${\rm Ai}(x)$ of the first kind which is the unique solution with $\lim_{x\rightarrow\infty}{\rm Ai}(x)=0$ of the differential equation ${\rm Ai}''(x)=x{\rm Ai}(x)$.
\end{thm}

Similar asymptotic results but for different combinatorial counting problems were obtained by Elvey Price et al. in \cite{EPFaWa}, where the (unusual) term $\exp\{cn^{\alpha}\}$ with $c$ some constant and $\alpha<1$ was called a {\it stretched exponential}. In fact, the method from \cite{EPFaWa} will also play a crucial role in the proof of our result.

As a consequence of our method of proof, we can also get an asymptotic result for the number of tree-child networks were all nodes (except the root) are bijectively labeled. We denote this number by $\widehat{TC}_{N}$ with $N$ the number of non-root nodes. Then, we have the following result which we formulate as corollary.

\begin{cor}\label{cor-lab}
The number of tree-child networks with $N$ non-root nodes which are all labeled satisfies
\[
\widehat{TC}_{N}=\left(\frac{3}{e^5}+o(1)\right)^{N/4}N^{5N/4},
\]
where $N$ runs through all odd positive integers.
\end{cor}

Recall that if $N$ is even, then a tree-child network with $N$ non-root nodes does not exist; see \cite{DiSeWe}.

\begin{Rem}
In \cite{DiSeWe}, the number of tree-child networks with all non-root nodes labeled was also considered and the authors showed that $N^{5N/4}$ is the dominating term in the main asymptotic growth term. Our above result shows that $(3/e^5)^{1/4}$ is the base of the exponential growth rate.
\end{Rem}

We next give a short sketch of the proof of Theorem~\ref{main-result}. First, let $\mathcal{TC}_{n,k}$ denote the set of tree-child networks with $n$ leaves and $k$ reticulation nodes and denote its cardinality by ${\rm TC}_{n,k}$. It is easy to see that $0\leq k\leq n-1$. For fixed values of $k$, the asymptotics of ${\rm TC}_{n,k}$ were derived in \cite{FuGiMa}:
\[
{\rm TC}_{n,k}\sim c_k\left(\frac{2}{e}\right)^nn^{n+2k-1},
\]
where $c_k>0$ is a computable constant.

However, in order to understand the asymptotics of ${\rm TC}_{n}$, it turns out that ${\rm TC}_{n,n-1}$ plays the most crucial role. More precisely, we will first show that
\begin{equation}\label{step-1}
{\rm TC}_n=\Theta({\rm TC}_{n,n-1}).
\end{equation}
Next, we will observe that ${\rm TC}_{n,n-1}/n!$ is actually contained in the OEIS\footnote{http://www.oeis.org/} as entry A213863 (with a shift). The latter sequence, say $a_n$, is defined as the number of words with letters $\{\omega_1,\ldots,\omega_n\}$ where each letter can be used exactly three times and in each prefix of the words, the number of occurrences of letter $\omega_i$ is either zero or if it is non-zero, then the letter $\omega_i$ must occur at least as often as the letter $\omega_j$ for all $j>i$. We give a bijective proof that ${\rm TC}_{n,n-1}/n!$ is indeed equal to $a_{n-1}$. Finally, for $a_n$ we will be able to find a recurrence to which the method from \cite{EPFaWa} can be applied.

The above will be done in the next three sections of the paper. Then, in Section~\ref{proof-cor}, we will give the proof of Corollary~\ref{cor-lab}. In Section~\ref{count-1}, we will derive the asymptotics of the number of {\it 1-component tree-child networks} for which explicit formulas were recently given in \cite{Zh2}. Again, stretched exponentials will occur in both the leaf-labeled and node-labeled case (even three of them in the latter). We will conclude the paper in Section~\ref{con} with a summary and some open problems.

\section{Tree-Child Networks with a Maximal Number of Reticulation Nodes}

In this section, we prove (\ref{step-1}). First, recall the following relation between $n,k,$ and $t$ of a phylogenetic network with $n$ leaves, $k$ reticulation nodes, and $t$ tree nodes:
\begin{equation}\label{fund-rel}
n+k=t+1.
\end{equation}
This follows from the fact that the sum of the out-degrees of nodes is equal to the sum of in-degrees; see \cite{DiSeWe} for details. 

Next, we say that a tree node is {\it free} if each of its children is either a tree node or a leaf. Moreover, we call an edge to a child of a free tree node a {\it free edge}. Then, we have the following lemma.
\begin{lmm}
Every tree-child network in $\mathcal{TC}_{n,k}$ has $n-k-1$ free tree nodes and thus $2(n-k-1)$ free edges.
\end{lmm}
\pf From (\ref{fund-rel}), we have that a tree-child network from $\mathcal{TC}_{n,k}$ has $n+k-1$ tree nodes. The two parents of every reticulation node are not free and due to the tree-child property, different reticulation nodes have different parents. Thus, the number of tree nodes which are not free is $2k$ from which the result follows.\qed

Using the previous lemma, we can show the following result.

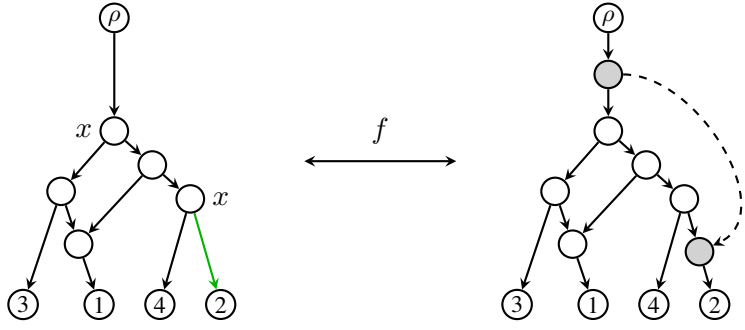
\begin{figure}
\begin{center}
\begin{tikzpicture}
\draw (0cm,0cm) node[inner sep=1.2pt,line width=0.8pt,draw,circle] (1) {{\footnotesize $\rho$}};
\draw (0cm,-1.5cm) node[minimum size=7.8pt,line width=0.8pt,draw,circle] (2) {};
\draw (-0.7cm,-2.3cm) node[minimum size=7.8pt,line width=0.8pt,draw,circle] (3) {};
\draw (-1.2cm,-3.8cm) node[inner sep=1.2pt,line width=0.8pt,draw,circle] (5) {{\footnotesize 3}};
\draw (-0.47cm,-3cm) node[minimum size=7.8pt,line width=0.8pt,draw,circle] (4) {};
\draw (-0.2cm,-3.8cm) node[inner sep=1.2pt,line width=0.8pt,draw,circle] (6) {{\footnotesize 1}};
\draw (0.5cm,-1.95cm) node[minimum size=7.8pt,line width=0.8pt,draw,circle] (7) {};
\draw (1cm,-2.4cm) node[minimum size=7.8pt,line width=0.8pt,draw,circle] (8) {};
\draw (0.6cm,-3.8cm) node[inner sep=1.2pt,line width=0.8pt,draw,circle] (9) {{\footnotesize 4}};
\draw (1.4cm,-3.8cm) node[inner sep=1.2pt,line width=0.8pt,draw,circle] (10) {{\footnotesize 2}};

\draw (-0.4cm,-1.5cm) node (11) {$x$};
\draw (1.4cm,-2.4cm) node (12) {$x$};

\draw[-stealth,line width=0.8pt] (1) -- (2);
\draw[-stealth,line width=0.8pt] (2) -- (3);
\draw[-stealth,line width=0.8pt] (3) -- (4);
\draw[-stealth,line width=0.8pt] (4) -- (6);
\draw[-stealth,line width=0.8pt] (3) -- (5);
\draw[-stealth,line width=0.8pt] (2) -- (7);
\draw[-stealth,line width=0.8pt] (7) -- (4);
\draw[-stealth,line width=0.8pt] (7) -- (8);
\draw[-stealth,line width=0.8pt] (8) -- (9);
\draw[-stealth,line width=0.8pt,color=black!30!green] (8) -- (10);

\draw[stealth-stealth,line width=0.8pt] (2.5cm,-1.9cm) -- (4.5cm,-1.9cm);
\draw (3.5cm,-1.5cm) node (0) {$f$};

\draw (6.5cm,0cm) node[inner sep=1.2pt,line width=0.8pt,draw,circle] (1) {{\footnotesize $\rho$}};
\draw (6.5cm,-0.75cm) node[minimum size=7.8pt,line width=0.8pt,draw,circle,fill=gray!35] (11) {};
\draw (6.5cm,-1.5cm) node[minimum size=7.8pt,line width=0.8pt,draw,circle] (2) {};
\draw (5.8cm,-2.3cm) node[minimum size=7.8pt,line width=0.8pt,draw,circle] (3) {};
\draw (5.3cm,-3.8cm) node[inner sep=1.2pt,line width=0.8pt,draw,circle] (5) {{\footnotesize 3}};
\draw (6.03cm,-3cm) node[minimum size=7.8pt,line width=0.8pt,draw,circle] (4) {};
\draw (6.3cm,-3.8cm) node[inner sep=1.2pt,line width=0.8pt,draw,circle] (6) {{\footnotesize 1}};
\draw (7cm,-1.95cm) node[minimum size=7.8pt,line width=0.8pt,draw,circle] (7) {};
\draw (7.5cm,-2.4cm) node[minimum size=7.8pt,line width=0.8pt,draw,circle] (8) {};
\draw (7.1cm,-3.8cm) node[inner sep=1.2pt,line width=0.8pt,draw,circle] (9) {{\footnotesize 4}};
\draw (7.7cm,-3.1cm) node[minimum size=7.8pt,line width=0.8pt,draw,circle,fill=gray!35] (12) {};
\draw (7.9cm,-3.8cm) node[inner sep=1.2pt,line width=0.8pt,draw,circle] (10) {{\footnotesize 2}};

\draw[-stealth,line width=0.8pt] (1) -- (11);
\draw[-stealth,line width=0.8pt] (11) -- (2);
\draw[-stealth,line width=0.8pt] (2) -- (3);
\draw[-stealth,line width=0.8pt] (3) -- (4);
\draw[-stealth,line width=0.8pt] (4) -- (6);
\draw[-stealth,line width=0.8pt] (3) -- (5);
\draw[-stealth,line width=0.8pt] (2) -- (7);
\draw[-stealth,line width=0.8pt] (7) -- (4);
\draw[-stealth,line width=0.8pt] (7) -- (8);
\draw[-stealth,line width=0.8pt] (8) -- (9);
\draw[-stealth,line width=0.8pt] (8) -- (12);
\draw[-stealth,line width=0.8pt] (12) -- (10);
\draw[dashed,-stealth,line width=0.8pt] (11) to [out=0,in=25] (12);
\end{tikzpicture}
\end{center}
\caption{The injective map $f$ from Lemma~\ref{growth}. (The two free nodes are marked by an $x$ and the picked free edge is in green.)}\label{inj-f}
\end{figure}

\begin{lmm}\label{growth}
For any $0\leq k\leq n-2$, we have $2(n-k-1){\rm TC}_{n,k}\leq{\rm TC}_{n,k+1}$.
\end{lmm}
\pf Define a map $f$ from a pair of a network from $\mathcal{TC}_{n,k}$ and a free edge of that network into $\mathcal{TC}_{n,k+1}$ as follows: insert a tree node between the root of the network and its child and add an edge from this node to a reticulation node which is inserted into the free edge; see Figure~\ref{inj-f}. Note that each image of $f$ is indeed a tree-child network from $\mathcal{TC}_{n,k+1}$ and that $f$ is injective. Moreover, because of the previous lemma, the range of $f$ has cardinality $2(n-k-1){\rm TC}_{n,k}$. From this the claim follows.\qed

We also need the following lemma.

\begin{figure}[t!]
\centering
\includegraphics[scale=1.0]{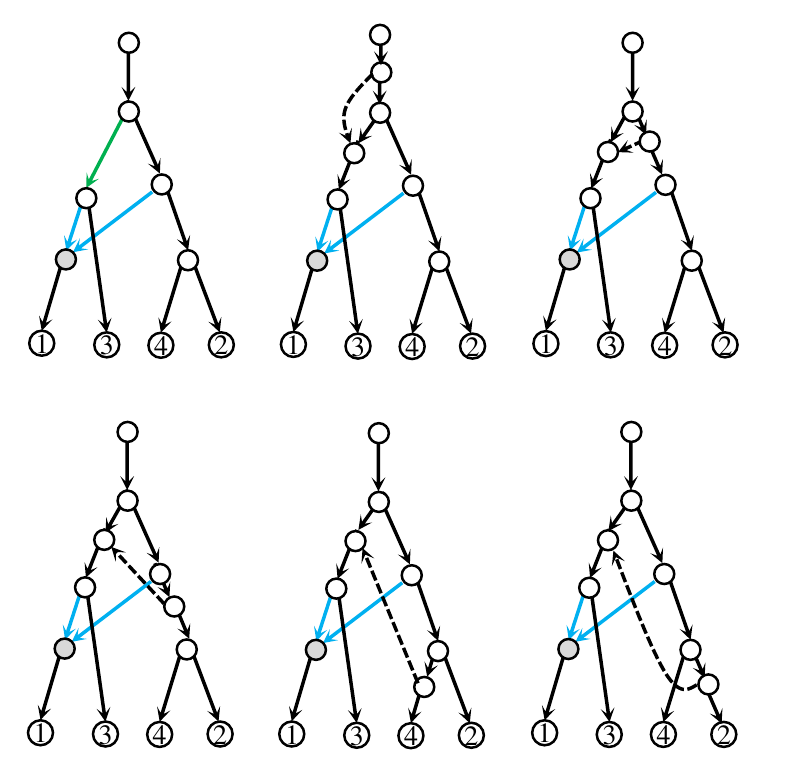}
\caption{Illustration of the proof of Lemma~\ref{lemma3} with $n=4$. In the top row, the leftmost tree-child network has $n-3$ ($=1$) reticulation nodes and hence contains $4$ free edges, one of which is colored green, and 7 other tree edges (black). Since the tree edges entering leaves 1 and 3 are below the green edge, we could not add an edge from either to form a DAG. Thus, we obtain  5 tree-child networks with $n-2$ reticulation nodes by adding an edge from each of the other 5 tree edges to the green free edge.  \label{expansion_fig}}
 \end{figure}

\begin{lmm}
\label{lemma3}
For $n\geq 3$, we have ${\rm TC}_{n, n-3}\geq \frac{1}{8}\times {\rm TC}_{n,n-2}$.
\end{lmm}
\pf First note that each network in $\mathcal{TC}_{n,n-3}$ contains $4$ free edges and $3n-4$ edges ending either in a tree node or leaf. (Such edges are subsequently called {\it tree edges}.) By picking a free edge and a (different) tree edge, then inserting a tree node $u$ into the picked tree edge and connecting it to a reticulation node $v$ which is inserted into the picked free edge, we obtain a network of ${\cal TC}_{n, n-2}$ if $u$ is not below $v$ (see Figure~\ref{expansion_fig}).  In this way, we get at most $4((3n-4)-1){\rm TC}_{n, n-3}$ networks with $n-2$ reticulation nodes, as we discard those that are not DAGs. 
On the other hand, each tree-child network of ${\cal TC}_{n, n-2}$ is obtained by this construction exactly $2(n-2)$ times.  Thus, $4(3n-5){\rm TC}_{n, n-3}\geq 2(n-2) {\rm TC}_{n, n-2}$. From this, since for $n\geq 3$,
\[
\frac{n-2}{3n-5}= \frac{1}{3}\left(1-\frac{1}{3n-5}\right)\geq \frac{1}{4},
\]
the inequality of the lemma follows.\qed

Now, we can prove the following quantitative version of (\ref{step-1}).

\begin{pro}\label{theta-bound-TCn}
For any $n\geq 3$, we have $\frac{25}{16}\times {\rm TC}_{n,n-1}\leq{\rm TC}_{n}\leq\sqrt{e}\times{\rm TC}_{n,n-1}$.
\end{pro}
\pf The lower bound follows from
\begin{equation}\label{id-tc-n}
{\rm TC}_{n}=\sum_{k=0}^{n-1}{\rm TC}_{n,k},
\end{equation}
the fact that ${\rm TC}_{n,n-1}=2{\rm TC}_{n,n-2}$ which was proved in \cite{Zh2}, and Lemma~\ref{lemma3}.

For the upper bound, by Lemma~\ref{growth} and iteration, we get
\begin{equation}\label{ub-TCnk}
{\rm TC}_{n,k}\leq\frac{1}{2^{n-1-k}(n-k-1)!}{\rm TC}_{n,n-1},\qquad (0\leq k\leq n-1).
\end{equation}
Thus, from (\ref{id-tc-n}),
\[
{\rm TC}_n\leq\left(\sum_{k=0}^{n-1}\frac{1}{2^{n-1-k}(n-k-1)!}\right){\rm TC}_{n,n-1}\leq\sqrt{e}{\rm TC}_{n,n-1}.
\]
This proves the claim.\qed

The above proposition reduces the problem of finding the main term of the asymptotics of ${\rm TC}_{n}$ to that of ${\rm TC}_{n,n-1}$ which is the number of tree-child networks with $n$ leaves and a maximal number of $n-1$ reticulation nodes.

These networks have a special structure which we discuss next. Recall that by Definition~\ref{def-tc}, for every node in the network there exists a path starting with that node and ending with a leaf whose intermediate nodes are all tree nodes. For tree-child networks with a maximal number of reticulation nodes, we have the following characterization.

\begin{lmm}
A tree-child network from $\mathcal{TC}_n$ has $n-1$ reticulation nodes if and only if the path from every node to a leaf whose intermediate nodes are all tree nodes is unique.
\end{lmm}
\pf First, assume that we have a tree-child network with $n$ leaves and $n-1$ reticulation nodes. Then, for different reticulation nodes, the paths from these nodes to leaves with all intermediate nodes being tree nodes end with different leaves. Moreover, the child of the root (which is a tree node) also has a path with all intermediate nodes being tree nodes which ends with a yet another leaf. Thus, we have already at least $n$ leaves and consequently, no node can have two paths with the claimed property because then the number of leaves would exceed $n$.

Next, assume that for every node there is a unique path to a leaf with all intermediate nodes being tree nodes. Consider first this path for the child of the root. Clearly, all intermediate nodes must be parents of reticulation nodes for otherwise an intermediate node would have two different paths to leaves with all intermediate nodes being tree nodes. Moreover, any reticulation node which is the child of an intermediate node on the path is followed by a tree node which again has a path to a leaf with all intermediate nodes being parents of reticulation nodes, etc. Clearly, this gives a network with $n$ leaves and exactly $n-1$ reticulation nodes.\qed

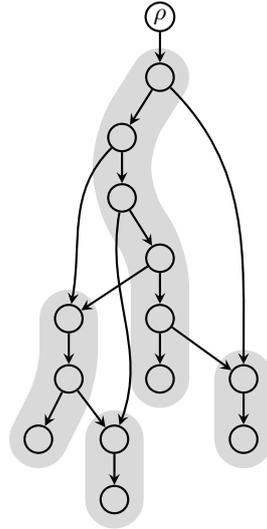
\begin{figure}
\begin{center}
\begin{tikzpicture}
\draw (0cm,0cm) node[inner sep=1.2pt,line width=0.8pt,draw,circle] (1) {{\footnotesize $\rho$}};
\draw (0cm,-0.8cm) node[minimum size=7.8pt,line width=0.8pt,draw,circle] (2) {};
\draw (-0.5cm,-1.6cm) node[minimum size=7.8pt,line width=0.8pt,draw,circle] (3) {};
\draw (-0.5cm,-2.4cm) node[minimum size=7.8pt,line width=0.8pt,draw,circle] (4) {};
\draw (0cm,-3.2cm) node[minimum size=7.8pt,line width=0.8pt,draw,circle] (5) {};
\draw (0cm,-4cm) node[minimum size=7.8pt,line width=0.8pt,draw,circle] (6) {};
\draw (0cm,-4.8cm) node[minimum size=7.8pt,line width=0.8pt,draw,circle] (7) {};
\draw (-1.2cm,-4cm) node[minimum size=7.8pt,line width=0.8pt,draw,circle] (8) {};
\draw (-1.2cm,-4.8cm) node[minimum size=7.8pt,line width=0.8pt,draw,circle] (9) {};
\draw (-1.6cm,-5.6cm) node[minimum size=7.8pt,line width=0.8pt,draw,circle] (10) {};
\draw (-0.6cm,-5.6cm) node[minimum size=7.8pt,line width=0.8pt,draw,circle] (11) {};
\draw (-0.6cm,-6.4cm) node[minimum size=7.8pt,line width=0.8pt,draw,circle] (12) {};
\draw (1.1cm,-4.8cm) node[minimum size=7.8pt,line width=0.8pt,draw,circle] (13) {};
\draw (1.1cm,-5.6cm) node[minimum size=7.8pt,line width=0.8pt,draw,circle] (14) {};

\draw[-stealth,line width=0.8pt] (1) -- (2);
\draw[-stealth,line width=0.8pt] (2) -- (3);
\draw[-stealth,line width=0.8pt] (3) -- (4);
\draw[-stealth,line width=0.8pt] (4) -- (5);
\draw[-stealth,line width=0.8pt] (5) -- (6);
\draw[-stealth,line width=0.8pt] (6) -- (7);
\draw[-stealth,line width=0.8pt] (3) to [out=-135,in=80] (8);
\draw[-stealth,line width=0.8pt] (5) -- (8);
\draw[-stealth,line width=0.8pt] (8) -- (9);
\draw[-stealth,line width=0.8pt] (9) -- (10);
\draw[-stealth,line width=0.8pt] (9) -- (11);
\draw[-stealth,line width=0.8pt] (4) to [out=-100,in=70] (11);
\draw[-stealth,line width=0.8pt] (11) -- (12);
\draw[-stealth,line width=0.8pt] (2) to [out=-45,in=90] (13);
\draw[-stealth,line width=0.8pt] (6) -- (13);
\draw[-stealth,line width=0.8pt] (13) -- (14);

\begin{pgfonlayer}{background}
\draw[rounded corners=1em,line width=2em,gray!30,cap=round] (2.center) -- (3.center) -- (4.center) -- (5.center) -- (6.center) -- (7.center);
\draw[rounded corners=1em,line width=2em,gray!30,cap=round] (8.center) -- (9.center) -- (10.center);
\draw[rounded corners=1em,line width=2em,gray!30,cap=round] (11.center) -- (12.center);
\draw[rounded corners=1em,line width=2em,gray!30,cap=round] (13.center) -- (14.center);
\end{pgfonlayer}
\end{tikzpicture}
\end{center}
\caption{A tree-child network with maximal number of reticulation nodes. (The path-components are highlighted and labels of leaves are removed; note that different labelings of the leaves lead to different tree-child networks.)}\label{max-retic}
\end{figure}

The second part of the proof of the above lemma explains the structure of tree-child networks with $n$ leaves and $n-1$ reticulation nodes: they consist of {\it path-components} starting with either a reticulation node or the child of the root. All other edges are between these components; see Figure~\ref{max-retic}. In particular, note that these networks have no symmetry and thus, if we remove the labels of the leaves, the number of all resulting (now unlabeled) networks is ${\rm TC}_{n,n-1}/n!$.

\section{A Class of Words and Recurrences for Their Counting Sequence}

We start with entry A213863 in the OEIS which is a counting sequence of certain words.
\begin{Def}\label{def-words}
Let ${\mathcal A}_n$ denote the class of words on the $n$ letters $\{\omega_1,\ldots,\omega_n\}$ such that in each word every letter occurs exactly $3$ times and in every prefix, the letter $\omega_i$ either has not yet occurred or if it has occurred, then the number of occurrences is at least as large as the number of occurrences of $\omega_j$ for all $j>i$. Moreover, denote by $a_n$ the cardinality of ${\mathcal A}_n$.
\end{Def}

In the OEIS, the first 16 terms of $a_n$ were given together with a brute-force Maple program to compute further terms (which becomes very slow beyond the $20$-th term). We recall the first $7$ terms: 
\[
\{a_n\}_{n\geq 1}=\{1, 7, 106, 2575, 87595, 3864040, 210455470, \ldots\}.
\]
In fact, it turns out that $a_{n-1}$ also counts the number of tree-child networks with $n$ leaves and $n-1$ reticulation nodes with the labels of leaves removed.

\begin{pro}\label{TC-and-an}
There is a bijection from the set of tree-child networks $\mathcal{TC}_{n,n-1}$ with labels removed to ${\mathcal A}_{n-1}$. Consequently, $a_{n-1}={\rm TC}_{n,n-1}/n!$.
\end{pro}
\pf We directly give the bijection. Therefore, start with an element from $\mathcal{TC}_{n,n-1}$ with labels removed; see Figure~\ref{max-retic} for an example. Recall that these objects are counted by ${\rm TC}_{n,n-1}/n!$ due to the lack of symmetry of the networks from $\mathcal{TC}_{n,n-1}$.

\begin{figure}
\begin{center}
\begin{tikzpicture}
\draw (0cm,0cm) node[inner sep=1.2pt,line width=0.8pt,draw,circle] (1) {{\footnotesize $\rho$}};
\draw (0cm,-0.8cm) node[minimum size=7.8pt,line width=0.8pt,draw,circle] (2) {};
\draw (-0.5cm,-1.6cm) node[minimum size=7.8pt,line width=0.8pt,draw,circle] (3) {};
\draw (-0.5cm,-2.4cm) node[minimum size=7.8pt,line width=0.8pt,draw,circle] (4) {};
\draw (0cm,-3.2cm) node[minimum size=7.8pt,line width=0.8pt,draw,circle] (5) {};
\draw (0cm,-4cm) node[minimum size=7.8pt,line width=0.8pt,draw,circle] (6) {};
\draw (0cm,-4.8cm) node[minimum size=7.8pt,line width=0.8pt,draw,circle] (7) {};
\draw (-1.2cm,-4cm) node[minimum size=7.8pt,line width=0.8pt,draw,circle] (8) {};
\draw (-1.2cm,-4.8cm) node[minimum size=7.8pt,line width=0.8pt,draw,circle] (9) {};
\draw (-1.6cm,-5.6cm) node[minimum size=7.8pt,line width=0.8pt,draw,circle] (10) {};
\draw (-0.6cm,-5.6cm) node[minimum size=7.8pt,line width=0.8pt,draw,circle] (11) {};
\draw (-0.6cm,-6.4cm) node[minimum size=7.8pt,line width=0.8pt,draw,circle] (12) {};
\draw (1.1cm,-4.8cm) node[minimum size=7.8pt,line width=0.8pt,draw,circle] (13) {};
\draw (1.1cm,-5.6cm) node[minimum size=7.8pt,line width=0.8pt,draw,circle] (14) {};

\draw[-stealth,line width=0.8pt] (1) -- (2);
\draw[-stealth,line width=0.8pt] (2) -- (3);
\draw[-stealth,line width=0.8pt] (3) -- (4);
\draw[-stealth,line width=0.8pt] (4) -- (5);
\draw[-stealth,line width=0.8pt] (5) -- (6);
\draw[-stealth,line width=0.8pt] (6) -- (7);
\draw[-stealth,line width=0.8pt] (3) to [out=-135,in=80] (8);
\draw[-stealth,line width=0.8pt] (5) -- (8);
\draw[-stealth,line width=0.8pt] (8) -- (9);
\draw[-stealth,line width=0.8pt] (9) -- (10);
\draw[-stealth,line width=0.8pt] (9) -- (11);
\draw[-stealth,line width=0.8pt] (4) to [out=-100,in=70] (11);
\draw[-stealth,line width=0.8pt] (11) -- (12);
\draw[-stealth,line width=0.8pt] (2) to [out=-45,in=90] (13);
\draw[-stealth,line width=0.8pt] (6) -- (13);
\draw[-stealth,line width=0.8pt] (13) -- (14);

\begin{pgfonlayer}{background}
\draw[rounded corners=1em,line width=2em,gray!30,cap=round] (2.center) -- (3.center) -- (4.center) -- (5.center) -- (6.center) -- (7.center);
\draw[rounded corners=1em,line width=2em,gray!30,cap=round] (8.center) -- (9.center) -- (10.center);
\draw[rounded corners=1em,line width=2em,gray!30,cap=round] (11.center) -- (12.center);
\draw[rounded corners=1em,line width=2em,gray!30,cap=round] (13.center) -- (14.center);
\end{pgfonlayer}

\draw (-1.1cm,-1.6cm) node (15) {$0$};
\draw (-1.9cm,-4.8cm) node (16) {$1$};
\draw (0cm,-6cm) node (17) {$3$};
\draw (1.7cm,-5.2cm) node (18) {$2$};

\draw[stealth-stealth,line width=0.8pt] (2.5cm,-3.2cm) -- (3.5cm,-3.2cm);

\draw (6cm,0cm) node[inner sep=1.2pt,line width=0.8pt,draw,circle] (1) {{\footnotesize $\rho$}};
\draw (6cm,-0.8cm) node[inner sep=1.2pt,line width=0.8pt,draw,circle] (2) {{\footnotesize $2$}};
\draw (5.5cm,-1.6cm) node[inner sep=1.2pt,line width=0.8pt,draw,circle] (3) {{\footnotesize $1$}};
\draw (5.5cm,-2.4cm) node[inner sep=1.2pt,line width=0.8pt,draw,circle] (4) {{\footnotesize $3$}};
\draw (6cm,-3.2cm) node[inner sep=1.2pt,line width=0.8pt,draw,circle] (5) {{\footnotesize $1$}};
\draw (6cm,-4cm) node[inner sep=1.2pt,line width=0.8pt,draw,circle] (6) {{\footnotesize $2$}};
\draw (6cm,-4.8cm) node[minimum size=7.8pt,line width=0.8pt,draw,circle] (7) {};
\draw (4.8cm,-4cm) node[inner sep=1.2pt,line width=0.8pt,draw,circle] (8) {{\footnotesize $1$}};
\draw (4.8cm,-4.8cm) node[inner sep=1.2pt,line width=0.8pt,draw,circle] (9) {{\footnotesize $3$}};
\draw (4.4cm,-5.6cm) node[minimum size=7.8pt,line width=0.8pt,draw,circle] (10) {};
\draw (5.4cm,-5.6cm) node[inner sep=1.2pt,line width=0.8pt,draw,circle] (11) {{\footnotesize $3$}};
\draw (5.4cm,-6.4cm) node[minimum size=7.8pt,line width=0.8pt,draw,circle] (12) {};
\draw (7.1cm,-4.8cm) node[inner sep=1.2pt,line width=0.8pt,draw,circle] (13) {{\footnotesize $2$}};
\draw (7.1cm,-5.6cm) node[minimum size=7.8pt,line width=0.8pt,draw,circle] (14) {};

\draw[-stealth,line width=0.8pt] (1) -- (2);
\draw[-stealth,line width=0.8pt] (2) -- (3);
\draw[-stealth,line width=0.8pt] (3) -- (4);
\draw[-stealth,line width=0.8pt] (4) -- (5);
\draw[-stealth,line width=0.8pt] (5) -- (6);
\draw[-stealth,line width=0.8pt] (6) -- (7);
\draw[-stealth,line width=0.8pt] (3) to [out=-135,in=80] (8);
\draw[-stealth,line width=0.8pt] (5) -- (8);
\draw[-stealth,line width=0.8pt] (8) -- (9);
\draw[-stealth,line width=0.8pt] (9) -- (10);
\draw[-stealth,line width=0.8pt] (9) -- (11);
\draw[-stealth,line width=0.8pt] (4) to [out=-100,in=70] (11);
\draw[-stealth,line width=0.8pt] (11) -- (12);
\draw[-stealth,line width=0.8pt] (2) to [out=-45,in=90] (13);
\draw[-stealth,line width=0.8pt] (6) -- (13);
\draw[-stealth,line width=0.8pt] (13) -- (14);

\begin{pgfonlayer}{background}
\draw[rounded corners=1em,line width=2em,gray!30,cap=round] (2.center) -- (3.center) -- (4.center) -- (5.center) -- (6.center) -- (7.center);
\draw[rounded corners=1em,line width=2em,gray!30,cap=round] (8.center) -- (9.center) -- (10.center);
\draw[rounded corners=1em,line width=2em,gray!30,cap=round] (11.center) -- (12.center);
\draw[rounded corners=1em,line width=2em,gray!30,cap=round] (13.center) -- (14.center);
\end{pgfonlayer}

\draw[stealth-stealth,line width=0.8pt] (8.5cm,-3.2cm) -- (9.5cm,-3.2cm);

\draw (11.3cm,-3.2cm) node (0) {$21312\vert13\vert2\vert3$};
\end{tikzpicture}
\end{center}
\caption{The bijection between unlabeled elements from $\mathcal{TC}_{n,n-1}$ and words from ${\mathcal A}_{n-1}$ from Proposition~\ref{TC-and-an}. (Here, $n=4$.)}\label{tc-to-words}
\end{figure}
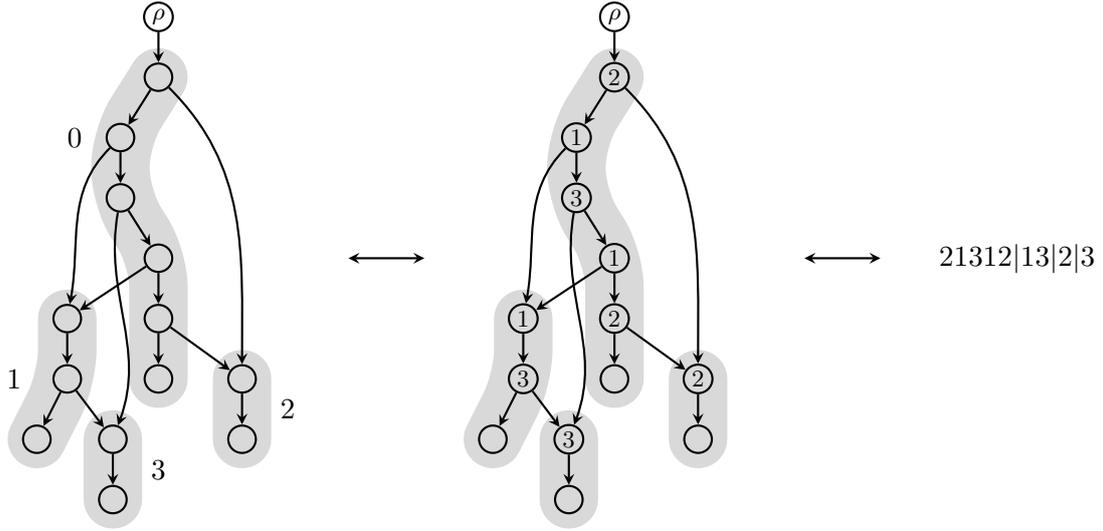

In the first step, we order the path-components of the chosen tree-child network. We do this inductively. First, the path-component of the child of the root receives index $0$. Assume that $k$ path-components have been indexed. Now, consider all un-indexed path-components whose first node (which is a reticulation node) has its two parents already in indexed path-components. If both parents are in the same path-component, then one is the descendant of the other; call that one the {\it second parent}; if both parents are in different path-components, then the parent in the path-component with the higher index is the {\it second parent}. Now order all the above chosen un-indexed path-components according to the indices of the path-components where the second parents are located and in case indices coincide, the ancestor relationship within the path-component of their second parents. Continue this until all path-components are indexed which will eventually happen because our networks are assumed to be connected; see the first part of Figure~\ref{tc-to-words}.

Now, we label the first node of every path-component of index $k>0$ together with its two parents by $k$; see the second part of Figure~\ref{tc-to-words}.

Finally, we read the labels of each path-component starting with the $0$-th one until we reach the last one; see the third part of Figure~\ref{tc-to-words}, where the separation line separates the strings from different path-components.

The resulting word is a word with letters $\{1,\ldots,n-1\}$ with each letter repeated exactly three times. Moreover, if a letter of the resulting word when read from the left occurs for the first time, then due to the above construction, no larger letter can have occurred already (at least) twice. Likewise, if a letter occurs for the second time, again no larger letter can have occurred already three times. Thus, the resulting words satisfies the property from Definition~\ref{def-words}.

Finally, it is straightforward to see that the above construction can be reversed. Thus, the resulting map is a bijection.\qed

For the sequence $a_n$, we can give now a recurrence.

\begin{pro}
Let $(b_{n,m})_{n,m\geq 1}$ be defined recursively as
\begin{equation}\label{rec-1-bnk}
b_{n,m}=(2n+m-2)\sum_{k=1}^{m}b_{n-1,k},\qquad (n\geq 2, 1\leq m\leq n)
\end{equation}
with initial conditions $b_{n,m}=0$ for $n<m$ and $b_{1,1}=1$. Then,
\[
a_n=\sum_{m\geq 1}b_{n,m}.
\]
\end{pro}
\pf First, note that any word in ${\mathcal A}_n$ must end with the letter $\omega_n$. Moreover, if one considers the suffix of a word from the second occurrence of $\omega_n$ to the last occurrence of $\omega_n$, then this suffix must have the form $\omega_n\omega_{n-i}\cdots\omega_{n-1}\omega_n$ for some $0\leq i\leq n-1$ (where here and in the sequel, for $i=0$ the suffix is understood to be $\omega_n\omega_n$).

Define now ${\mathcal A}_{n,i}$ as the set of all words in ${\mathcal A}_n$ with suffix $\omega_n\omega_{n-i}\cdots\omega_{n-1}\omega_n$ for some $0\leq i\leq n-1$. Then,
\[
{\mathcal A}_n={\mathcal A}_{n,0}\cup\cdots\cup{\mathcal A}_{n,n-1}
\]
and this union is disjoint. Let $a_{n,i}$ be the cardinality of ${\mathcal A}_{n,i}$.

First, consider ${\mathcal A}_{n,0}$. Here, we can generate all words from words in ${\mathcal A}_{n-1}$ by placing $\omega_n\omega_n$ at the end and inserting a third $\omega_n$ anywhere into the word. This gives,
\begin{equation}\label{rec-part-1}
a_{n,0}=(3n-2)a_{n-1},
\end{equation}
since the third $\omega_n$ has exactly $3n-2$ positions to choose from.

Next, consider ${\mathcal A}_{n,i}$ with $1\leq i\leq n-1$. By definition, all words in this set have the suffix $\omega_n\omega_{n-i}\cdots\omega_{n-1}\omega_n$. Moreover, if we remove the three $\omega_n$'s, then the resulting word ends with $\omega_{n-i}\cdots\omega_{n-1}$ and is thus contained in
\[
{\mathcal A}_{n-1,i-1}\cup\cdots\cup{\mathcal A}_{n-1,n-2}.
\]
Conversely, if we insert three $\omega_n$'s into a word of the above set by placing one at the beginning and one at the end of the suffix $\omega_{n-i}\cdots\omega_{n-1}$ and inserting the third $\omega_n$ anywhere before this, then we create all words from ${\mathcal A}_{n,i}$. This shows that
\begin{equation}\label{rec-part-2}
a_{n,i}=(3n-2-i)\sum_{k=i-1}^{n-2}a_{n-1,k}
\end{equation}
because the third $\omega_n$ has exactly $3n-2-i$ positions to choose from.

Finally, by setting
\[
b_{n,i}:=a_{n,n-i}
\]
we obtain the claimed recurrence from (\ref{rec-part-1}) and (\ref{rec-part-2}) with the claimed initial conditions (which are easily verified).\qed

Using this recurrence, the first, e.g., $1000$ terms of $a_n$ can be computed with Maple in a few seconds.

The recurrence for $b_{n,m}$ can also be written in the following equivalent form.
\begin{cor}\label{rec-2-bnk}
We have,
\[
b_{n,m}=\frac{2n+m-2}{2n+m-3}b_{n,m-1}+(2n+m-2)b_{n-1,m},\qquad (n\geq 2, 0\leq m\leq n)
\]
with initial conditions $b_{n,m}=0$ for $m>n$, $b_{n,-1}=0$ for $n\geq 1$ and $b_{1,m}=0$ for $m\geq -1$ except for $m=1$ where $b_{1,1}=1$.
\end{cor}
\pf We have,
\[
\frac{b_{n,m}}{2n+m-2}=\sum_{k=1}^{m}b_{n-1,k}.
\]
Thus,
\[
\frac{b_{n,m}}{2n+m-2}-\frac{b_{n,m-1}}{2n+m-3}=b_{n-1,m}
\]
which can be re-arranged to the claimed recurrence. Finally, the initial conditions are easily adjusted to the current recurrence.\qed

\section{Asymptotic Analysis and Proof of Theorem~\ref{main-result}}\label{sec-proof-mt}

In this section, we will do the asymptotic analysis of $a_n$. As a warm-up, we first start with an upper bound which follows from (\ref{rec-1-bnk}).

\begin{lmm}\label{crude-bound-an}
We have,
\[
a_n={\mathcal O}\left(n^{-1}\left(\frac{12}{e}\right)^nn^{n}\right).
\]
\end{lmm}
\pf Define a sequence $g_{n,m}$ recursively by
\[
g_{n,m}=\sum_{k=1}^{m}g_{n-1,k},\qquad(n\geq 2,1\leq m\leq n)
\]
with initial conditions $g_{n,m}=0$ for $n<m$ and $g_{1,1}=1$. Then, by induction and using (\ref{rec-1-bnk}), we get that
\begin{equation}\label{upp-bound}
3^nn!g_{n,m}\geq b_{n,m},\qquad (n,m\geq 1).
\end{equation}

Next, note that $g_{n,m}$ are the ballot numbers which are found in the OEIS as entry A009766 and admit the closed form expression
\[
g_{n,m}=\frac{n-m+1}{n}\binom{n+m-2}{n-1}
\]
Thus,
\[
\sum_{m\geq 1}g_{n,m}=g_{n+1,n+1}=\frac{1}{n+1}\binom{2n}{n}=C_n,
\]
where $C_n$ denotes the $n$-th Catalan number, which by Stirling's formula has the asymptotics
\[
C_n\sim\frac{4^{n}}{\sqrt{\pi n^3}}.
\]

Finally, from (\ref{upp-bound}), we obtain that
\[
a_n=\sum_{m\geq 1}b_{n,m}\leq 3^nn!\sum_{m\geq 1}g_{n,m}=3^nn!C_n\sim\frac{\sqrt{2}}{n}\left(\frac{12}{e}\right)^n n^{n}.
\]
This gives the claimed result.\qed

Note that by using (\ref{step-1}) and Proposition~\ref{TC-and-an}, this upper bound implies that
\[
{\rm TC}_n={\mathcal O}\left(n^{-3/2}\left(\frac{12}{e^2}\right)^nn^{2n}\right).
\]
This already captures the right exponential growth order; compare with Theorem~\ref{main-result}. However, we are still missing the stretched exponential (and the polynomial term) which requires a much more subtle method. We will explain this method next.

The method we are going to use was introduced in \cite{EPFaWa}, where it was first applied (amongst other things) to another problem on DAGs, namely, the counting of relaxed binary trees. (Despite our problem being also a counting problem on DAGs, it came as a surprise for us that the the same method applies to both problems.) Indeed, the heuristic for the method from \cite{EPFaWa} comes from lattice path theory and is well-explained in Section~3 in \cite{EPFaWa}.

The main reason why the method from \cite{EPFaWa} applies in our context is that our recurrence from Corollary~\ref{rec-2-bnk} has the same shape as the one for relaxed binary trees which was given in Proposition~2.6 in \cite{EPFaWa}. The only difference is that the coefficients on the right-hand side are slightly different: in the recurrence for relaxed binary trees in \cite{EPFaWa}, the coefficients are $1$ and $m+1$, whereas our coefficients are $(2n+m-2)/(2n+m-3)$ and $2n+m-2$. But note that we are interested in the asymptotics of $b_{n,n}$ since by (\ref{rec-1-bnk}) and initial conditions:
\begin{equation}\label{rel-b-a}
b_{n,n}=(3n-2)a_{n-1}.
\end{equation}
Now, for $m\approx n$, our coefficients become $(2n+m-2)/(2n+m-3)\approx 1$ and $2n+m-2\approx 3n$ compared to the above coefficients for the recurrence of relaxed binary trees which become $1$ and $m+1\approx n$. Thus, we have good reason to believe that the method from \cite{EPFaWa} will apply to our sequence when divided by $3^n$. This will turn out to be indeed be the case.

We will now use the steps from \cite{EPFaWa}. The first step was to define
\[
d_{n,m}:=\frac{1}{3^{(n+m)/2}\left(\frac{n+m}{2}\right)!}b_{(n+m)/2,(n-m)/2},\qquad (n\geq 4, 0\leq m\leq n, n-m\ \text{even})
\]
and $d_{n,m}:=0$ if $n-m$ is odd. Note that $b_{n,n}=3^nn!d_{2n,0}$. Moreover, from Corollary~\ref{rec-2-bnk},
\begin{equation}\label{rec-dnm}
d_{n,m}=\frac{3n+m-4}{3n+m-6}d_{n-1,m+1}+\frac{3n+m-4}{3n+3m}d_{n-1,m-1},\qquad (n\geq 3, m\geq 0)
\end{equation}
with initial conditions $d_{n,-1}=0$ for $n\geq 2$, $d_{2,m}=0$ for $m\geq 1$ and $d_{2,0}=1/3$.

The next step is to assume that
\[
d_{n,m}\approx h(n)f\left((m+1)n^{-1/3}\right),
\]
where $h(n)\approx 2^n\mu^{n^{1/3}}$ with a (so far unknown) constant $\mu$. We repeat some of the heuristic arguments from Section~3 in \cite{EPFaWa} which show how to choose $\mu$ and $f(\cdot)$. The reason why we do this (instead of just referring to \cite{EPFaWa}) is because these arguments (i) show how the Airy function shows up and (ii) make it easier to understand the two propositions below. First, set $s(n)=h(n)/h(n-1)$ which satisfies
\begin{equation}\label{exp-s}
s(n)=2+\frac{2\log\mu}{3}n^{-2/3}+{\mathcal O}(n^{-1}).
\end{equation}
Now, plugging everything into (\ref{rec-dnm}), we obtain that
\begin{align}
s(n)f\left((m+1)n^{-1/3}\right)=\frac{3n+m-4}{3n+m-6}&f\left((m+2)(n-1)^{-1/3}\right)\nonumber\\
&\qquad+\frac{3n+m-4}{3n+3m}f\left(m(n-1)^{-1/3}\right).\label{eq-f}
\end{align}
Pick $m=\kappa n^{1/3}-1$. Then,  $
\frac{3n+m-4}{3n+m-6}=1+{\mathcal O}(n^{-1})$ and  $\frac{3n+m-4}{3n+3m}=1-\frac{2}{3}\kappa n^{-2/3}+{\mathcal O}(n^{-1})$.
Since
\[
f\left((m+2)(n-1)^{-1/3}\right)=f(\kappa)+f'(\kappa)n^{-1/3}+f''(\kappa)\frac{n^{-2/3}}{2}+{\mathcal O}(n^{-1})
\]
and
\[
f\left(m(n-1)^{-1/3}\right)=f(\kappa)-f'(\kappa)n^{-1/3}+f''(\kappa)\frac{n^{-2/3}}{2}+{\mathcal O}(n^{-1}),
\]
we obtain by plugging the above into (\ref{eq-f}) and rearranging
\[
\left(f''(\kappa)-\frac{2}{3}(\log\mu+\kappa)f(\kappa)\right)n^{-2/3}+{\mathcal O}(n^{-1})=0.
\]
Thus,
\[
f''(\kappa)=\frac{2}{3}(\log\mu+\kappa)f(\kappa)
\]
which under some regularity conditions on $f(\cdot)$ has the solution
\[
f(\kappa)=b{\rm Ai}\left(c(\log\mu+\kappa)\right),
\]
where $b$ is a suitable constant and $c=(2/3)^{1/3}$. If we further impose $f(0)=0$, then $c\log\mu=a_1$ and thus our function $f(\cdot)$ should be chosen as
\begin{equation}\label{shape-f}
f(\kappa)\approx{\rm Ai}\left(a_1+c\kappa\right).
\end{equation}

Now, the main idea of the approach in \cite{EPFaWa} is to choose $f(\cdot)$ of the form (\ref{shape-f}) with suitable multiplicative factors and suitable versions of (\ref{exp-s}) such that (\ref{eq-f}) holds as an inequality with both $\leq$ and $\geq$. Then, this gives lower and upper bounds for $d_{n,m}$ and thus $d_{2n,0}$ which then yields the desired upper and lower bound for $b_{n,n}$ and thus $a_n$ via (\ref{rel-b-a}).

We give the two results containing these inequalities next; compare with Lemma~4.2 and Lemma~4.4 in \cite{EPFaWa} for the corresponding results for relaxed binary trees.

\begin{pro}\label{lower-bound}
Define
\[
\tilde{X}_{n,m}:=\left(1-\frac{m^2}{3n}-\frac{25m}{18n}\right){\rm Ai}\left(a_1+c(m+1)n^{-1/3}\right),
\]
where $c=(2/3)^{1/3}$ and
\[
\tilde{s}(n):=2+\frac{2^{2/3}a_1}{3^{2/3}n^{2/3}}-\frac{2}{3n}-\frac{1}{n^{7/6}}.
\]
Then, for fixed $\epsilon>0$ and $n$ large enough, we have
\begin{equation}\label{lb}
\tilde{s}(n)\tilde{X}_{n,m}\leq\frac{3n+m-4}{3n+m-6}\tilde{X}_{n-1,m+1}+\frac{3n+m-4}{3n+3m}\tilde{X}_{n-1,m-1}
\end{equation}
for all $0\leq m<n^{2/3-\epsilon}$.
\end{pro}

\begin{pro}\label{upper-bound}
Let $\eta>1/18$ and define
\[
\hat{X}_{n,m}:=\left(1-\frac{m^2}{3n}-\frac{25m}{18n}+\eta\frac{m^4}{n^2}\right){\rm Ai}\left(a_1+c(m+1)n^{-1/3}\right),
\]
where $c=(2/3)^{1/3}$ and
\[
\hat{s}(n):=2+\frac{2^{2/3}a_1}{3^{2/3}n^{2/3}}-\frac{2}{3n}+\frac{1}{n^{7/6}}.
\]
Then, for fixed $\epsilon>0$ and $n$ large enough, we have
\begin{equation}\label{ub}
\hat{s}(n)\hat{X}_{n,m}\geq\frac{3n+m-4}{3n+m-6}\hat{X}_{n-1,m+1}+\frac{3n+m-4}{3n+3m}\hat{X}_{n-1,m-1}
\end{equation}
for all $0\leq m<n^{1-\epsilon}$.
\end{pro}

\vspace*{0.3cm}
\noindent{\it Proof of Proposition~\ref{lower-bound} and Proposition~\ref{upper-bound}.} For the proof, we use the Maple worksheet which was cited in \cite{EPFaWa} (item [26] in the list of references) and which the authors made available online. We only have to modify it to our situation; see the first autor's personal webpage for the modified version\footnote{http://web.math.nctu.edu.tw/mfuchs}.

First, for the claim in Proposition~\ref{lower-bound}, we set
\[
P_{n,m}:=-\tilde{s}(n)\tilde{X}_{n,m}+\frac{3n+m-4}{3n+m-6}\tilde{X}_{n-1,m+1}+\frac{3n+m-4}{3n+3m}\tilde{X}_{n-1,m-1}.
\]
Using the above mentioned Maple worksheet, by series expansion and the differential equation of the Airy function, this double-sequence can be written as
\begin{align*}
P_{n,m}=&{\rm Ai}(\alpha)\left(\frac{1}{n^{7/6}}-\frac{2^{5/3}a_1m}{3^{8/3}n^{5/3}}-\frac{23m^2}{27n^2}-\frac{2^{5/3}\cdot5a_1m^3}{3^{11/3}n^{8/3}}
-\frac{37m^4}{81n^3}+\frac{449m^5}{3645n^4}+\cdots\right)+\\
&{\rm Ai}'(\alpha)\Bigg(\frac{2^{1/3}}{3^{1/3}n^{3/2}}-\frac{2813m}{2^{2/3}\cdot3^{10/3}\cdot5n^{7/3}}-\frac{53\cdot2^{4/3}m^2}{3^{10/3}n^{7/3}}\\
&\qquad\qquad-\frac{2^{7/3}m^3}{3^{7/3}n^{7/3}}
-\frac{2^{1/3}m^4}{3^{7/3}n^{10/3}}-\frac{163\cdot2^{1/3}m^5}{3^{16/3}\cdot5n^{13/3}}+\cdots\Bigg),
\end{align*}
where $\alpha:=a_1+c(m+1)n^{-1/3}$. This has exactly the same shape as the expansion in the proof of Lemma~4.2 in \cite{EPFaWa} with the sole exception that the second term in the bracket behind ${\rm Ai}'$ is different but since this term is (i) now slightly smaller and (ii) was shown to be asymptotically not relevant in the proof of Lemma~4.2 in \cite{EPFaWa}, this does not matter. Thus, we can use the arguments from Lemma~4.2 in \cite{EPFaWa} to show that $P_{n,m}\geq 0$ for all large enough $n$ and $0\leq m<n^{2/3-\epsilon}$.

Similarly, Proposition~\ref{upper-bound} is proved with the Maple worksheet and the arguments from the proof of Lemma~4.4 in \cite{EPFaWa}.\qed

Now, equipped with the above two propositions, we can prove Theorem~\ref{main-result}, where the first part of the proof again heavily borrows from Section~4 in \cite{EPFaWa}.

\vspace*{0.3cm}\noindent{\it Proof of Theorem~\ref{main-result}.} First, define $\tilde{h}(n)=\tilde{s}(n)\tilde{h}(n-1)$ with $\tilde{h}(1):=1$. Then, (\ref{lb}) becomes
\[
\tilde{h}(n)\tilde{X}_{n,m}\leq\frac{3n+m-4}{3n+m-6}\tilde{h}(n-1)\tilde{X}_{n-1,m+1}+\frac{3n+m-4}{3n+3m}\tilde{h}(n-1)\tilde{X}_{n-1,m-1}
\]
for $n$ large and $0\leq m<n^{2/3-\epsilon}$. Note also that for $n$ large, $\tilde{X}_{n,m}$ is negative for $n^{2/3-\epsilon}\leq m\leq n$. Using this and induction, we obtain that
\[
d_{n,m}=\Omega\left(\tilde{h}(n)\max\{\tilde{X}_{n,m},0\}\right)
\]
for $n$ large and $0\leq m\leq n$. Thus,
\begin{align*}
b_{n,n}=3^nn!d_{2n,0}&=\Omega\left(3^nn!\tilde{h}(2n){\rm Ai}\left(a_1+c2^{-1/3}n^{-1/3}\right)\right)\\
&=\Omega\left(3^nn!\left(\prod_{\ell=2}^{2n}\tilde{s}(\ell)\right){\rm Ai}\left(a_1+c2^{-1/3}n^{-1/3}\right)\right),
\end{align*}
where the last step follows by iteration. Now, straightforward asymptotic expansion yields
\begin{align*}
\prod_{\ell=2}^{2n}\tilde{s}(\ell)&=\prod_{\ell=2}^{2n}\left(2+\frac{2^{2/3}a_1}{3^{2/3}\ell^{2/3}}-\frac{2}{3\ell}-\frac{1}{\ell^{7/6}}\right)\\
&=4^n\exp\left\{\frac{a_1}{2^{1/3}\cdot3^{2/3}}\sum_{\ell=2}^{2n}\ell^{-2/3}-3^{-1}\sum_{\ell=2}^{2n}\ell^{-1}+{\mathcal O}(1)\right\}\\
&=\Omega\left(n^{-1/3}e^{a_1(3n)^{1/3}}4^n\right)
\end{align*}
and
\[
{\rm Ai}\left(a_1+c2^{-1/3}n^{-1/3}\right)=c2^{-1/3}{\rm Ai}'(a_1)n^{-1/3}+{\mathcal O}\left(n^{-2/3}\right).
\]
Consequently,
\begin{equation}\label{lb-bnn}
b_{n,n}=\Omega\left(n!n^{-2/3}e^{a_1(3n)^{1/3}}12^n\right).
\end{equation}

Next, we use Proposition~\ref{upper-bound} to show a matching upper bound. Therefore, we again define $\hat{h}(n)=\hat{s}(n)\hat{h}(n-1)$ with $\hat{h}(1)=1$. Then, from (\ref{ub}),
\[
\hat{h}(n)\hat{X}_{n,m}\geq\frac{3n+m-4}{3n+m-6}\hat{h}(n-1)\hat{X}_{n-1,m+1}+\frac{3n+m-4}{3n+3m}\hat{h}(n-1)\hat{X}_{n-1,m-1},
\]
for $n$ large enough and $0\leq m<n^{1-\epsilon}$. Now, define a new sequence $\hat{d}_{n,m}$ which satisfies (\ref{rec-dnm}) for $0\leq m<n^{1-\epsilon}$ and $\hat{d}_{n,m}:=0$ otherwise. Then, by induction
\[
\hat{d}_{n,m}={\mathcal O}\left(\hat{h}(n)\hat{X}_{n,m}\right)
\]
for $n$ large and $0\leq m\leq n$. Thus, with the same arguments as above, we obtain that
\[
\hat{d}_{2n,0}={\mathcal O}\left(n^{-2/3}e^{a_1(3n)^{1/3}}4^n\right).
\]
Next, using the similar arguments as in \cite{EPFaWa} (which were based on lattice paths theory), one can show that for large $n$, we have $d_{2n,0}\leq 2\hat{d}_{2n,0}$; see the Appendix for details. Thus,
\begin{equation}\label{ub-bnn}
b_{n,n}=n!3^nd_{2n,0}={\mathcal O}\left(n!n^{-2/3}e^{a_1(3n)^{1/3}}12^n\right).
\end{equation}

Combining (\ref{lb-bnn}) and (\ref{ub-bnn}) now yields
\[
b_{n,n}=\Theta\left(n!n^{-2/3}e^{a_1(3n)^{1/3}}12^n\right).
\]
Next, from (\ref{rel-b-a}),
\begin{equation}\label{bound-an}
a_{n-1}=\Theta\left(n!n^{-5/3}e^{a_1(3n)^{1/3}}12^n\right)
\end{equation}
and finally by Proposition~\ref{theta-bound-TCn} and Proposition~\ref{TC-and-an},
\[
{\rm TC}_n=\Theta\left(n!^2n^{-5/3}e^{a_1(3n)^{1/3}}12^n\right)
\]
from which the claimed result follows by Stirling's formula.\qed

\section{Proof of Corollary~\ref{cor-lab}}\label{proof-cor}

In this section, we consider tree-child networks with $N$ non-root nodes which are bijectively labeled, where $N$ is odd (see the remark below Corollary~\ref{cor-lab}). 

We start with the following relationship between their number $\widehat{{\rm TC}}_N$ and ${\rm TC}_{n,k}$.

\begin{lmm}\label{formula-node-lab}
We have,
\begin{equation}\label{formula-cN}
\widehat{{\rm TC}}_N=\sum_{n=\lceil (N+3)/4\rceil}^{(N+1)/2}\frac{N!}{n!}\times{\rm TC}_{n,(N+1)/2-n}.
\end{equation}
\end{lmm}

\pf This is proved with the following two results from \cite{DiSeWe}.
\begin{itemize}
\item[(i)] The number $n$ of leaves, the number $k$ of reticulation nodes, and the number $N$ of non-root nodes satisfy:
\[
n+k=\frac{N+1}{2};
\]
see equation (5) in \cite{DiSeWe}.
\item[(ii)] The descendant sets of any two non-root nodes in a tree-child network are distinct; see Lemma 5.1 in \cite{DiSeWe}. (Here, the descendant set of a node is the set of leaves which are descendant from the node.)
\end{itemize}

More precisely, it follows from (ii) that in order to generate all tree-child networks with $N$ labeled non-root nodes, we can take a tree-child network with $n$ labeled leaves and $N$ non-root nodes, choose $n$ of the $N$ labels and use them to relabel the leaves such that the order between labels is preserved, and finally distribute the remaining $N-n$ labels arbitrarily to the non-leaf and non-root nodes. Here, because of (i), $n$ ranges from $\lceil(N+3)/4\rceil$ to $(N+1)/2$. Moreover, again by (i), the number of tree-child networks with $n$ labeled leaves and $N$ non-root nodes is given by ${\rm TC}_{n,(N+1)/2-n}$.

This shows the claimed result.\qed

Corollary~\ref{cor-lab} is now a consequence of the following two propositions.
\begin{pro}\label{lb-cN}
We have,
\[
\widehat{{\rm TC}}_N=\Omega\left(N^{1/12}e^{a_1(3N/4)^{1/3}}\left(\frac{3}{e^5}\right)^{N/4}
N^{5N/4}\right).
\]
\end{pro}
\pf Set
\[
\tilde{n}=\left\lceil\frac{N+3}{4}\right\rceil.
\]
Clearly, we have
\[
\widehat{{\rm TC}}_N\geq \frac{N!}{\tilde{n}!}{\rm TC}_{\tilde{n},(N+1)/2-\tilde{n}}.
\]
Note that
\[
{\rm TC}_{\tilde{n},(N+1)/2-\tilde{n}}=\begin{cases}{\rm TC}_{\tilde{n},\tilde{n}-1},&\text{if}\ N\equiv 1\ \text{mod}\ 4;\\
{\rm TC}_{\tilde{n},\tilde{n}-1}/2,&\text{if}\ N\equiv 3\ \text{mod}\ 4,\end{cases}
\]
where the second case follows from  $2{\rm TC}_{n,n-2}={\rm TC}_{n,n-1}$ which was proved in \cite{Zh2}. Thus, by Proposition~\ref{TC-and-an},
\[
\widehat{{\rm TC}}_N=\Omega\left(N!a_{\tilde{n}-1}\right).
\]
Now using (\ref{bound-an}), Stirling's formula and straightforward computation gives the claim.\qed

\begin{pro}
We have,
\[
\widehat{{\rm TC}}_N={\mathcal O}\left(N^2e^{\sqrt{3N}/2}\left(\frac{3}{e^5}\right)^{N/4}N^{5N/4}\right).
\]
\end{pro}
\pf By plugging (\ref{ub-TCnk}) into (\ref{formula-cN}), we have
\begin{align*}
\widehat{{\rm TC}}_N&\leq N!\sum_{n=\lceil(N+3)/4\rceil}^{(N+1)/2}\frac{{\rm TC}_{n,n-1}}{n!2^{2n-1-(N+1)/2}(2n-1-(N+1)/2)!}\\
&=N!2^{(N+3)/2}\sum_{n=\lceil(N+3)/4\rceil}^{(N+1)/2}\frac{a_{n-1}}{4^n(2n-1-(N+1)/2)!},
\end{align*}
where the last step follows from Proposition~\ref{TC-and-an}. Now, from Lemma~\ref{crude-bound-an}, we have
\begin{equation}\label{cru-bound-an}
a_{n-1}={\mathcal O}(n!n^{-5/2}12^n)
\end{equation}
and plugging this into the estimate for $\widehat{{\rm TC}}_{N}$ above, we get
\[
\widehat{{\rm TC}}_N={\mathcal O}\left(N!2^{N/2}N^{-5/2}\sum_{n=\lceil(N+3)/4\rceil}^{(N+1)/2}\frac{n!3^{n}}{(2n-1-(N+1)/2)!}\right).
\]
The terms inside the sum are maximal at
\[
\tilde{n}=\left\lfloor\frac{N}{4}+\frac{11}{8}+\frac{\sqrt{12N+61}}{8}\right\rfloor.
\]
Thus,
\[
\widehat{{\rm TC}}_N={\mathcal O}\left(N!2^{N/2}N^{-3/2}\frac{\tilde{n}!3^{\tilde{n}}}{(2\tilde{n}-1-(N+1)/2)!}\right)
\]
from which the claim follows by straightforward computation.\qed

\begin{Rem} When using (\ref{bound-an}) instead of (\ref{cru-bound-an}) and applying the Laplace method (see the next section), the upper bound in the last lemma can be improved. However, the subexponential large gap between the bounds in the last two propositions still remains even with such an improved upper bound.
\end{Rem}

\section{Asymptotic Count of 1-Component Tree-Child Networks}\label{count-1}

In this section, we consider 1-component tree-child networks which were defined in \cite{Zh2} as follows.

\begin{Def} A tree-child network is called 1-component tree-child network if every reticulation node is directly followed by a leaf.
\end{Def}

Denote by $1$-${\rm TC}_{n,k}$ the number of 1-component tree-child networks with $n$ leaves and $k$ reticulation nodes and by $1$-${\rm TC}_n$ the number of all 1-component tree-child networks with $n$ leaves. These number have easy, explicit formulas as was shown in Section 4.1 in \cite{Zh2}.
\begin{thm}[\cite{Zh2}]\label{for-1-TC}
The number of 1-component tree-child networks with $n$ leaves and $k$ reticulation nodes is given by
\begin{equation}\label{1TCnk}
\text{$1$-${\rm TC}_{n,k}$}=\binom{n}{k}\frac{(2n-2)!}{2^{n-1}(n-k-1)!},\qquad (0\leq k\leq n-1).
\end{equation}
Consequently,
\[
\text{$1$-${\rm TC}_n$}=\sum_{k=0}^{n-1}\binom{n}{k}\frac{(2n-2)!}{2^{n-1}(n-k-1)!}.
\]
\end{thm}

From the above formula for $1$-${\rm TC}_n$, we can get the first order asymptotics for this number with the Laplace method (which is the real-analytic version of the saddle point method); see, e.g., Chapter 9 in Graham et al. \cite{GrKnPa}.

\begin{thm}
The number of 1-component tree-child networks with $n$ leaves satisfies
\[
\text{$1$-${\rm TC}_n$}\sim\frac{1}{2\sqrt{e}}n^{-5/4}e^{2\sqrt{n}}\left(\frac{2}{e^2}\right)^nn^{2n}.
\]
\end{thm}

\begin{Rem}
Comparing with the result for the total number of tree-child networks with $n$ leaves from Theorem~\ref{main-result}, we see that the main term $n^{2n}$ is the same but the above exponential growth term is much smaller. Moreover, note there is again a stretched exponential in the asymptotics.
\end{Rem}

\pf By Theorem~\ref{for-1-TC},
\[
\text{$1$-${\rm TC}_n$}=\frac{n!(2n-2)!}{2^{n-1}}\sum_{k=0}^{n-1}\frac{1}{k!(n-k)!(n-k-1)!}.
\]
Stirling's formula gives
\begin{equation}\label{asymp-fac}
\frac{n!(2n-2)!}{2^{n-1}}\sim\sqrt{2}\pi n^{-1}\left(\frac{2}{e^{3}}\right)^nn^{3n}.
\end{equation}
Thus, we only need the asymptotics of
\[
S:=\sum_{k=0}^{n-1}\frac{1}{k!(n-k)!(n-k-1)!}
\]
which follows by a standard application of the Laplace method.

First observe that $1/(k!(n-k)!(n-k-1)!)$ is increasing for $k\leq n-\sqrt{n+1}$ and decreasing for $k>n-\sqrt{n+1}$. Moreover, by straightforward expansion,
\[
\frac{1}{k!(n-k)!(n-k-1)!}=\frac{1}{2\pi\sqrt{2e\pi}}n^{-1/2}e^{2\sqrt{n}}\left(\frac{e}{n}\right)^ne^{-x^2/\sqrt{n}}\left(1+{\mathcal O}\left(\frac{1+\vert x\vert}{\sqrt{n}}\right)\right),
\]
uniformly for $\vert x\vert\leq n^{3/8}$ where $k=n-\sqrt{n}+x$. Thus,
\begin{align*}
S&\sim\frac{1}{2\pi\sqrt{2e\pi}}n^{-1/2}e^{2\sqrt{n}}\left(\frac{e}{n}\right)^n\sum_{-n^{3/8}\leq x\leq n^{3/8}}e^{-x^2/\sqrt{n}}\left(1+{\mathcal O}\left(\frac{1+\vert x\vert}{\sqrt{n}}\right)\right)\\
&\sim\frac{1}{2\pi\sqrt{2e\pi}}n^{-1/2}e^{2\sqrt{n}}\left(\frac{e}{n}\right)^n\int_{-\infty}^{\infty}e^{-x^2/\sqrt{n}}\left(1+{\mathcal O}\left(\frac{1+\vert x\vert}{\sqrt{n}}\right)\right){\rm d}x\\
&\sim\frac{1}{2\pi\sqrt{2e}}n^{-1/4}e^{2\sqrt{n}}\left(\frac{e}{n}\right)^n\left(1+{\mathcal O}\left(\frac{1}{\sqrt{n}}\right)\right),
\end{align*}
where we used here standard arguments to justify tail-pruning and attachment; see \cite{GrKnPa}.

Combining this with (\ref{asymp-fac}) gives the claimed result.\qed

Next, we consider 1-component tree-child networks with all non-root nodes bijectively labeled. Denote by $1$-$\widehat{{\rm TC}}_{N}$ their number, where $N$ is the number of non-root nodes. Then, with the same proof method as in  Lemma~\ref{formula-node-lab}, we have
\begin{equation}\label{1hatTC}
1\text{-}\widehat{{\rm TC}}_N=\sum_{n=\lceil(N+3)/4\rceil}^{(N+1)/2}\binom{N}{n}(N-n)!\times\left(1\text{-}{\rm TC}_{n,(N+1)/2-n}\right).
\end{equation}
From this, by another application of the Laplace method, we obtain the following result.

\begin{thm}
The number of 1-component tree-child networks with all non-root nodes labeled satisfies
\[
1\text{-}\widehat{{\rm TC}}_N\sim2^{9/8}e^{-1/(32)}N^{-7/8}e^{(2N)^{3/4}/2+\sqrt{2N}/(16)-3(2N)^{1/4}/(64)}\left(\frac{1}{2e^5}\right)^{N/4}N^{5N/4}.
\]
\end{thm}

\begin{Rem}
Again the main term $N^{5N/4}$ is the same as in Corollary~\ref{cor-lab} but the exponential growth term is different. Also, note that now there are three stretched exponentials in the asymptotics.
\end{Rem}

\pf Plugging (\ref{1TCnk}) into (\ref{1hatTC}), we obtain that
\[
1\text{-}\widehat{TC}_N=2N!\sum_{n=\lceil(N+3)/4\rceil}^{(N+1)/2}\frac{(2n-2)!}{((N+1)/2-n)!(2n-(N+1)/2)!(2n-(N+3)/2)!2^n}.
\]
Thus, the asymptotics will follow from that of
\[
S:=\sum_{n=\lceil(N+3)/4\rceil}^{(N+1)/2}\frac{(2n-2)!}{((N+1)/2-n)!(2n-(N+1)/2)!(2n-(N+3)/2)!2^n}
\]
which will be deduced by another application of the Laplace method.

First, observe that the terms of $S$ become maximal at $\tilde{n}+1$ where $\tilde{n}$ is the largest integer such that
\[
\frac{8(N+1-2\tilde{n})(2\tilde{n}-1)\tilde{n}}{(N-4\tilde{n}+1)(N-4\tilde{n}-1)^2(N-4\tilde{n}-3)}\leq 1.
\]
From this, by bootstrapping,
\[
\tilde{n}=\frac{N}{4}+\frac{(2N)^{3/4}}{8}+\frac{\sqrt{2N}}{32}-\frac{9(2N)^{1/4}}{256}+{\mathcal O}(1).
\]
Now, by a long computation (which is best done with the help of, e.g., Maple):
\begin{align*}
&\frac{(2n-2)!}{((N+1)/2-n)!(2n-(N+1)/2)!(2n-(N+3)/2)!2^n}\\
&\qquad=2^{-3/8}\pi^{-1}e^{-1/(32)}N^{-7/4}q(n)\left(\frac{1}{2e}\right)^{N/4}N^{N/4}e^{-x^2N^{-3/4}}\left(1+{\mathcal O}\left(\frac{1}{N^{1/4}}+\frac{ x^2}{N}\right)\right)
\end{align*}
uniformly for $\vert x\vert\leq n^{2/5}$, where
\[
k=\frac{N}{4}+\frac{(2N)^{3/4}}{8}+\frac{\sqrt{2N}}{32}-\frac{9(2N)^{1/4}}{256}+x
\]
and
\[
q(n):=e^{(2N)^{3/4}/2+\sqrt{2N}/(16)-3(2N)^{1/4}/(64)}.
\]

Thus,
\begin{align*}
S&\sim\frac{q(n)}{2^{3/8}\pi e^{1/(32)}}N^{-7/4}\left(\frac{1}{2e}\right)^{N/4}N^{N/4}\int_{-\infty}^{\infty}e^{-x^2N^{-3/4}}\left(1+{\mathcal O}\left(\frac{1}{N^{1/4}}+\frac{x^2}{N}\right)\right){\rm d}x\\
&\sim\frac{q(n)}{2^{3/8}\sqrt{\pi} e^{1/(32)}}N^{-11/8}\left(\frac{1}{2e}\right)^{N/4}N^{N/4},
\end{align*}
where we again used standard arguments.

The claimed result follows now from the above asymptotic of $S$ multiplied with the asymptotic of $2N!$.\qed

\section{Conclusion}\label{con}

In \cite{DiSeWe}, the authors proved that $n^{2n}$ resp.~$N^{5N/4}$ are the dominating terms in the main asymptotic growth term of the number of tree-child networks with $n$ labeled leaves resp. $N$ labeled non-root nodes. Moreover, they asked to find in both cases the exponential growth terms.

In this paper, we answered these questions. Moreover, in the leaf-labeled case, our result gives all terms of the main asymptotic growth term up to a constant. Interestingly, this growth term contains a stretched exponential. Furthermore, we showed that stretched exponentials are also present in the asymptotics of the number of 1-component tree-child networks for both the leaf-labeled and node-labeled case. These numbers were considerably easier to analyze due to explicit formulas from \cite{Zh2}. Using these formulas, we even obtained the first-order asymptotics (and the method used in this paper would be capable of giving further terms in the asymptotic expansion, too).

Several interesting questions remain open. First, can the result of Theorem~\ref{main-result} be improved to a first-order asymptotic result, e.g., is there a constant $\gamma$ such that
\begin{align}\label{ref-asymp}
{\rm TC}_n\sim\gamma n^{-2/3}e^{a_1(3n)^{1/3}}\left(\frac{12}{e^2}\right)^{n}n^{2n}?
\end{align}
In fact, this was also discussed by the authors from \cite{EPFaWa} for the combinatorial counting problems in their paper. However, as pointed out in Section~3.4 in \cite{EPFaWa}, the approach from that paper (which played also a crucial role in the proof of Theorem~\ref{main-result}) is incapable of providing such a refined result. Thus, one needs to come up with a new approach to be able to prove such  results. However, note that such an improved approach, in our situation would ``only" give the first-order asymptotics of $a_n$, which is still not enough for a proof of (\ref{ref-asymp}) because one in addition would also need to improve Proposition~\ref{theta-bound-TCn} to a first-order asymptotic result. (Note that the bound in Proposition~\ref{theta-bound-TCn} is already quite tight since the constant in the upper bound is $\sqrt{e}=1.64872\ldots$ compared to the constant $25/16=1.5625$ in the lower bound.)

A second interesting question is whether the result from Corollary~\ref{cor-lab} for tree-child networks with all non-root nodes labeled can be improved? Is it possible to find further terms of the main asymptotic growth term? Also, are again stretched exponential(s) present?

A final interesting question is whether similar asymptotic results as those from Theorem~\ref{main-result} and Corollary~\ref{cor-lab} can also be proved for other classes of phylogenetic networks, in particular for normal networks? Here, a tree-child network is normal if for each reticulation node neither of its parents is a child of the other. In fact, in \cite{DiSeWe}, the authors showed that the number of normal networks in the leaf-labeled resp. node-labeled case has also the terms $n^{2n}$ resp.~$N^{5n/4}$ in their main asymptotic growth terms. Moreover, they also showed  that the number of normal networks in both cases is a small-o of the number of tree-child networks. Again, they asked for the exponential growth rates (in particular whether or not they coincide with the ones for tree-child networks) and we can now go one step further and ask whether there are also stretched exponentials in the main asymptotic growth term of the number of normal networks in the leaf-labeled and node-labeled case? Note that our method, rather surprisingly, does not work for the number of normal networks since a result similar to Proposition~\ref{theta-bound-TCn} does not seem to hold. In fact, the data from \cite{Zh1} suggests that for fixed $n$, the number of normal networks with $n$ leaves and $k$ reticulation nodes is unimodal in contrast to tree-child networks where this sequence grows at least exponentially; compare with Lemma~\ref{growth}.

\section*{Acknowledgments}

This work was started when the first author visited the third author at the National University of Singapore in November 2019. He thanks him for hospitality and financial support. The second author joined the research when he visited the Institute of Statistical Sciences, Academia Sinica in December 2019. He thanks Hsien-Kuei Hwang for hospitality and financial support. 

We also thank both referees for a careful reading and many helpful remarks/suggestions leading to a substantial improvement of the presentation. In particular, the short proof of Lemma~\ref{crude-bound-an} was suggested by one of the referees.

\section*{Appendix}

The goal of this appendix is to prove the claimed inequality $d_{2n,0}\leq2\hat{d}_{2n,0}$ which was stated at the end of Section~\ref{sec-proof-mt}. As mentioned there, the proof will use lattice paths theory and proceeds as in Section~4.2 in \cite{EPFaWa}. In order to avoid repetition, we will only point out differences.

First, from (\ref{rec-dnm}), we see that $d_{n,m}$ is the number of weighted lattice path starting at $(0,0)$ and ending at $(n,m)$ with up-steps of the form $(1,1)$ starting at $(a,b)$ having weight $(3a+b)/(3a+3b+6)$ and down-steps of the form $(1,-1)$ also starting at $(a,b)$ having weight $(3a+b-2)/(3a+b-4)$. Each lattice path is counted with its weight which is the product of all weights of its up-steps and down-steps.

Define $p_{\ell,m,2n}$ to be the number of lattice paths of the above type which start at $(\ell,m)$ and end at $(2n,0)$. Note that this number satisfies the recurrence
\begin{equation}\label{rec-plmn}
p_{\ell,m,2n}=\frac{3\ell+m}{3\ell+3m+6}p_{\ell+1,m+1,2n}+\frac{3\ell+m-2}{3\ell+m-4}p_{\ell+1,m-1,2n}\qquad (\ell,m\geq 0)
\end{equation}
with initial conditions $p_{\ell,-1,2n}=0$ and $p_{2n,m,2n}=0$ for $m\geq 1$ and $p_{2n,0,2n}=-1$.

A crucial step in the proof of the corresponding inequality to the one above in Section~4.2 of \cite{EPFaWa} was Lemma~4.5 which in our situation reads as follows.

\begin{lmm}
We have,
\[
\frac{p_{\ell,j,2n}}{(j+1)^2}\geq\frac{p_{\ell,k,2n}}{(k+1)^2},
\]
where $0\leq j<k\leq\ell\leq 2n$ and $2\vert k-j$.
\end{lmm}
\pf First note that it is sufficient to prove that
\begin{equation}\label{claimed-in}
\frac{p_{\ell,m-1,2n}}{m^2}-\frac{p_{\ell,m+1,2n}}{(m+2)^2}\geq 0
\end{equation}
for all $1\leq m\leq\ell-1\leq 2n-1$. We do this by reverse induction on $\ell$ and using (\ref{rec-plmn}). 

For $\ell=2n$ the claim is trivial. So, assume that it holds for $\ell+1$ and all $m$. We are going to prove it for $\ell$ and all $m$. 

In order to do so, we plug (\ref{rec-plmn}) into (\ref{claimed-in}) which gives
\begin{align*}
L:=&\frac{(3\ell+m-1)p_{\ell+1,m,2n}}{(3\ell+3m+3)m^2}+\frac{(3\ell+m-3)p_{\ell+1,m-2,2n}}{(3\ell+m-5)m^2}\\
&\qquad\quad-\frac{(3\ell+m+1)p_{\ell+1,m+2,2n}}{(3\ell+3m+9)(m+2)^2}-\frac{(3\ell+m-1)p_{\ell+1,m,2n}}{(3\ell+m-3)(m+2)^2}.
\end{align*}
Now, from the induction hypothesis, we have
\[
p_{\ell+1,m-2,2n}\geq\frac{(m-1)^2p_{\ell+1,m,2n}}{(m+1)^2}\qquad\text{and}\qquad p_{\ell+1,m+2,2n}\leq\frac{(m+3)^2p_{\ell+1,m,2n}}{(m+1)^2}.
\]
Plugging this into $L$ and re-arranging the obtained expression gives
\[
L=\frac{4(p_1(m)\ell^4+p_2(m)\ell^3+p_3(m)\ell^2+p_4(m)\ell+p_5(m))p_{\ell+1,m,2n}}{3(\ell+m+3)(\ell+m+1)(3\ell+m-3)(3\ell+m-5)(m+2)^2(m+1)^2m^2},
\]
where
\begin{align*}
p_1(m)&=54m+54;\\
p_2(m)&=9m^4+36m^3+117m^2+72m+54;\\
p_3(m)&=6m^5-24m^4-204m^3-474m^2-414m-282;\\
p_4(m)&=m^6-16m^5-36m^4+180m^3+537m^2+380m+138\\
p_5(m)&=20m^5+84m^4+24m^3-184m^2-108m+36.
\end{align*}
By computing all roots of these five polynomials (for instance with Maple), we see that $p_i(m)\geq 0$ for $1\leq i\leq 5$ and all $m\geq 18$ and consequently, $L\geq 0$ for $18\leq m\leq \ell-1$. As for the remaining $17$ cases of $m$, again it can be checked by Maple that
\[
p_1(m)\ell^4+p_2(m)\ell^3+p_3(m)\ell^2+p_4(m)\ell+p_5(m)\geq 0
\]
for $m=1,\ldots,17$ and $\ell\geq m+1$. 

Overall, we have $L\geq 0$ for $1\leq m\leq\ell-1$ which is the desired result.\qed

The rest of the proof of the inequality $d_{n,m}\leq 2\hat{d}_{n,m}$ proceeds now exactly as in Section~4.2 in \cite{EPFaWa} with the only difference that the bound 
\[
d_{2n,2m}\leq\binom{2n}{n+m},
\]
which was used in the proof in Section~4.2 of \cite{EPFaWa} where it followed from the fact that all weights are $\leq 1$, holds here despite the weight of down-steps being slightly larger than $1$. Indeed, by combining corresponding up-steps and down-steps, it is easy to see that the multiplied weight of such pairs of steps is $\leq 1$. We leave the (easy) computation to the reader.
\end{document}